\documentclass[a4paper,12pt]{article}

\def\BibTeX{{\rmfamily B\kern-.05em{\scshape i\kern-.025em b}\kern-.08em \TeX}}
\usepackage{verbatim}
\usepackage{psfig}
\usepackage[dvips]{graphicx}
\usepackage{amsfonts}
\usepackage{verbatim}
\newcommand{\cov}[2]{
  \mbox{\textsf{cov}}\left(#1,#2\right)
}
\newcommand{\var}[1]{
  \mbox{\textsf{var}}\left(#1\right)
}

\newtheorem{theorem}{\textbf{Theorem}}

\newtheorem{corollary}{\textbf{Corollary}}
\newtheorem{remark}{{\textbf{Remark}}}

\newcommand{\Proof}[1]{
{\noindent
\textit{Proof:} #1 \hfill $\Box$
}
}

\title{Bounds for covariances and variances of truncated random variables} 

\author{
N.~Hemachandra 
\thanks{IE and OR Interdisciplinary Programme, IIT Bombay; 
email: \texttt{nh@iitb.ac.in} }\\
Indian Institute of Technology Bombay, Mumbai, 400 076, India \\
and \\
V.~Cheriyan 
\thanks{e-mail: \texttt{vcheriyan@dc.com}. Most of this work was done when VC
was a Research Associate supported by Project 96IR001 at IITB. }\\
Deloitte Consulting LLP, Irving, TX - 75038, USA.}
\date{\today}

\begin{document}
\maketitle

\begin{abstract}
We show that a lower bound for covariance of $\min(X_1,X_2)$ and $\max(X_1,X_2)$
 is $\cov{X_1}{X_2}$ and an upper bound for variance of \\
$\min(X_2,\max(X,X_1))$ 
is $\var{X} + \var{X_1} +\var{X_2}$ generalizing previous results. 
We also characterize the cases where these  bounds are sharp.
\end{abstract}


{\em MSC 60 }

\begin{center}
\large{
\underline{October 2001. Revised version: October 2002}

}
\vspace{1cm}
Technical report 02\_2001, IE and OR Interdisciplinary Programme, IIT Bombay, Mumbai, 400 076, India. 
\end{center}
\normalsize

\newpage
\section{Introduction}
\label{secn:intro}

Truncation is known to be a basic tool of probability; see, for example, \cite{Chow88}. Chow and Teicher also identify two ways of truncating a random 
variable $X$: 
(a) $Y := \min(c,\max(X,a))$ 
(b) $Y^\prime := X1_{\{ a \leq X \leq c\}}$ 
for $a$ and $c$ constants such that $-\infty \leq a < c \leq \infty$. In this 
section, we first 
give two situations where truncated random variables arise while analyzing 
models. Then, we summarize the known results about the bounds of variances 
and covariances of such truncated random variables. 

Let $\{D_n\}_1^N, N \le \infty$ be a discrete time $Z^+$ valued process with 
$D_n$ being the demand for an item in period $n$. 
Let $I_0$ be the initial inventory and and for $n \ge 1$, let $I_{n}$ be the 
inventory at the end of the period $n$. Let $y_n$ be 
the inventory on hand at the beginning of the period $n$, 
 being possibly supplemented by ordering some quantity. 
The excess inventory in 
each period is charged a holding cost of $h$ per unit time per item. There is a 
penalty cost of $p$ for each item of demand that is not met; such portions of 
demand are back-ordered, {\it i.e.}, met in subsequent periods when 
items are available. 
 The unit purchase price of item is $c$. 
Then, $c_n$ the total cost in period $n$, involves truncated random variables: 
\[ c_n = c(y_n - I_{n-1}) + h \max(D_n - y_n, 0) + p \max(y_n - D_n, 0) \]

\noindent
The planning horizon $N$ could be finite or infinite and optimal policies that 
minimize different costs like 
total expected discounted cost or long run average cost are sought. 
A related model is to assume that unfulfilled orders are not back-ordered but  
are lost. Then the amount supplied can be represented as a truncated random 
variable. 
Work in the area of finding the optimal policies in such situations 
is largely initiated by Scarf \cite{Scarf} and Iglehart \cite{Iglehart}; 
these and many of the later developments are summarized in 
handbooks like \cite{Graves} and books like \cite{Bertsekas}. 

Our second example is from queueing models: Let $A_n$ be the time between the 
arrival times of customers $n$ and $n+1$ to a single server queue, 
and $S_n$ be the service time of customer $n$, $n \ge 1$
 If customers are served according to first-come-first-serve discipline, 
({\it i.e.,} in the order of their arrival), then waiting time of 
$n^{th}$ customer $W_n$, which is the time between $n^{th}$ customer's arrival 
and departure after service, is related by Lindley's recursion \cite{Lindley}: 
$W_{n+1} = \max (W_n + S_n - A_n, 0)$ with $W_0 = 0$. If $\{A_n\}_{n \ge 1}$ 
and $\{S_n\}_{n \ge 1}$ are independent {\it i.i.d.} sequences, 
one can find the asymptotic behaviour of the queue by analyzing the associated
 random walk \cite{Lindley}; see also, Asmussen~\cite{Asmussen}, 
Wolff~\cite{Wolff}, {\it etc.}, for details and related work. 

Chow and Studden \cite{Chow69} (using the notation given in the beginning), 
 note that for finite reals, $a$ and $c$, $\var{Y} \leq \var{X}$ 
while no comparable relationship exits between $\var{Y^\prime}$ and $\var{X}$. 
In fact, they show, among other things, that
\[ E[(Y-E[Y|\mathcal{G}])^2|\mathcal{G}] \leq E[(X-E[X|\mathcal{G}])^2|\mathcal{G}] \mbox{\it\ a.s.} \]
where, now $Y := \max(a,\min(X,b))$,  $\mathcal{G}$ is a sub $\sigma$-field and $a$ and $b$ are random variables that are $\mathcal{G}$ measurable, and $X,Y$ are integrable.

 
Consider an example: Let $\Omega = \{\omega_1, \omega_2, \omega_3 \}$ with 
$P(\omega_1) = \frac{1}{2}$ and $P(\omega_2) = P(\omega_3) = \frac{1}{4}$. 
Let $X$ be a random variable on $\Omega$ such that $X(\omega_1) = 0$ and 
$X(\omega_2) = X(\omega_3) = 1$; $X_2(\omega_1) = 2, X_2(\omega_2) = 0$ and 
$X_2(\omega_3) = 1$; and $X_1 = 2$ on $\Omega$. 
Let $X(X_1,X_2) := \min(X_2, \max(X,X_1))$. 
We have $\var{X} = \frac{1}{4}$ and $\var{X(X_1,X_2)}= \frac{11}{16}$, so that the variance of $X(X_1,X_2)$ is greater than the variance of $X$. 
We  give a sharp bound for variances of these type of random variables; 
it turns out that this bound also bounds variances of other truncated random 
variables obtained from $X,X_1$ and $X_2$. 


Let $X_1$ and $X_2$ be two random variables and let $Y := \min(X_1,X_2)$ and $Z
:= \max(X_1,X_2)$. 
One can view these $Y$ and $Z$ as order statistics of $X_1$ 
and $X_2$. In the context of finding estimators for dependent random variables, 
Papadatos \cite{Papadatos} has shown that 
$\cov{Y}{Z} \ge \cov{X_1}{X_2}$ if $X_1$ and $X_2$ have same law and are 
possibly dependent. 
We show that this result is also true even if distributions of  
$X_1$ and $X_2$ are different and  also characterize the cases when this bound 
is sharp. We use this result later to have upper bounds for variances of 
truncated random variables.  

\section{A lower bound for covariances}
\label{covsec}
Recall that for random variables $X_1$ and $X_2$ we define $Y := \min(X_1,X_2)$
 and $Z := \max(X_1,X_2)$. 
\begin{theorem}
\label{cov-bound}
 If $\var{X_i} < \infty,\; i = 1, 2$, then\\
\begin{equation}
\cov{Y}{Z} \geq \cov{X_1}{X_2} 
\label{cov-rel}
\end{equation}
with equality iff either $X_1 \ge X_2~a.s.$ or $X_2 \ge X_1~a.s.$
\end{theorem}

\Proof{ 
Since both $\{\min(X_1,X_2)\}^2$ and $\{\max(X_1,X_2)\}^2$ are upper bounded by
 $X_1^2 +X_2^2$  pointwise, both $Y$ and $Z$ have finite second moments. 

Next, we have, pointwise, 
\begin{equation}
\label{sumeqn}
Y+Z = X_1 + X_2
\end{equation}
\begin{equation}
\label{prodeqn}
Y Z = X_1 X_2
\end{equation}
Then,
\begin{eqnarray}
\cov{Y}{Z} - \cov{X_1}{X_2} &=& E[YZ]-E[Y]E[Z]-E[X_1X_2]+E[X_1]E[X_2] \nonumber \\
{} &=& E[X_1]E[X_2] - E[Y]E[Z] \label{covdiff}
\end{eqnarray}
from (\ref{prodeqn}), since $E[X_1X_2]<\infty$. Since $E[Y]E[X_1]$ is finite, 
we can add and subtract it to (\ref{covdiff}) to have,  
\begin{eqnarray}
\cov{Y}{Z} - \cov{X_1}{X_2} &=& E[X_1]E[X_2]-E[Y]E[Z]-E[X_1]E[Y]
\nonumber \\
{} &{}&  + E[X_1]E[Y] \nonumber\\
{} &=& E[X_1]\{E[X_2]-E[Y]\}-E[Y]\{E[Z]-E[X_1]\} \nonumber \\
\label{covcoveqn}
{} &=& ( E[X_1]-E[Y])(E[X_2]-E[Y])
\end{eqnarray}
because $E[Z]-E[X_1] = E[X_2]-E[Y]$ from (\ref{sumeqn}). As $X_i \geq Y, i =1,2$,
both the terms on the RHS of (\ref{covcoveqn}) are positive and thus 
(\ref{cov-bound})  follows. 
 
Next, equality exits in (\ref{cov-rel}) {\em iff} at least one term of 
(\ref{covcoveqn}) is zero. 
If $X_1\geq X_2$ {\it\ a.s.} ($X_2\geq X_1$ {\it\ a.s.}) then $Y=X_2$ ($Y=X_1$)
and hence equality holds.
On the other hand, if $E[X_1]-E[Y]=0$, {\em i.e.}, $\int (X_1-Y) dP =0$,
 then $X_1 = Y$ {\it\ a.s.} as $X_1 \geq Y$, so that $X_2 \geq X_1$. 
Similarly, if $E[X_2]-E[Y] =0$, then $X_1 \geq X_2~a.s$.
}

\begin{remark}
 From this, it follows that $Y$ and $Z$ are positively correlated
if $X_1$ and $X_2$ are so; in particular if $X_1$ and $X_2$ are independent,
then $Y$ and $Z$ are positively correlated. 
\end{remark}

\section{Bounds for variances}
\label{varsec}
We now consider the relationship between variances of random variables 
and those obtained by truncating them.
\begin{theorem}
\label{var-bound1}
 $\var{Y} + \var{Z} \leq \var{X_1} + \var{X_2}$ 
This is valid even if $\var{X_i}, i = 1,2$ is infinite for some 
$i = 1,2$.
\end{theorem}
\Proof{
{\em Case (i)} Suppose $\var{X_i} < \infty, i = 1,2$.
Then from (\ref{sumeqn})
\begin{equation}
\label{covcoveqn2}
\var{Y} + \var{Z} + 2(\cov{Y}{Z} -\cov{X_1}{X_2}) = \var{X_1} + \var{X_2}
\end{equation}
From Theorem (\ref{cov-bound}), $\cov{Y}{Z}-\cov{X_1}{X_2}$ is positive and 
hence, the result follows. 
\\
{\em Case (ii)} If variance of some $X_i$'s, $i=1,2$ is not finite, then
variances of $Y$ and $Z$ could be infinite and in that sense the inequality
holds. Note LHS is well defined.
}

We now characterize the cases when equality holds in the above.
 From (\ref{covcoveqn2}) and Theorem \ref{cov-bound} above, 
we have,
\begin{theorem}
\label{var-eqt} 
Suppose variance of $X_i, i=1,2$ is finite. Then,
\begin{equation}
\label{equaleqn}
\var{Y} + \var{Z} = \var{X_1} + \var{X_2}
\end{equation}
iff either $X_1\geq X_2$ {\it\ a.s.} or $X_2\geq X_1$ {\it\ a.s.}
\end{theorem}

\begin{remark}
If $X_1 \ge X_2$ {\it a.s.} or $X_2 \ge X_1$ {\it a.s.} then, (\ref{equaleqn}) 
holds. However, if (\ref{equaleqn}) holds, in general, it does not follow that
either $X_1\geq X_2$ {\it\ a.s.} or $X_2\geq X_1$ {\it\ a.s.} Consider $X$ with $P(X=i) =
\frac{1}{2ci^3}, i = \cdots,-3,-2,2,3,\cdots$ where $c=\sum_{i\geq 2}
\frac{1}{i^3}$. $X$ does not have finite variance. If $X_2 = d$, for a real
number $d$, then $Y$ and $Z$ also does not have finite variance. 
(\ref{equaleqn}) holds in the sense that both sides of it are infinite but we 
do not have $X \geq d$ {\it a.s.} or $d \geq X$ {\it a.s.}
\end{remark}

We would like to find conditions when variance of the truncated random variable
equals sum of variances of the random variables involved. 

\begin{corollary}
\label{cor1}
 Suppose variance of $X_i, i = 1,2$, is finite. Then 
\begin{equation}
\label{eqn9}
\var{Y} = \var{X_1} + \var{X_2}
\end{equation}
iff $Z$ is a constant {\it a.s.}, say $c$, and either
\begin{equation}
\label{condeqn}
c \geq X_1 \mbox{\it\ a.s,} \mbox{ or } c \geq X_2 \mbox{\it\ a.s.}
\end{equation}
\end{corollary}

\Proof{
If (\ref{eqn9}) holds, then using Theorem \ref{var-bound1} we have 
$\var{Z} = 0$,
{\it i.e.}, $Z$ is a constant {\it a.s.}, say $c$. From Theorem 
\ref{var-eqt} 
we have $X_2 \geq X_1$ {\it a.s.}, so that $X_2 = c$,
{\it a.s.}, or $X_1 \geq X_2$ {\it a.s.} so that $X_1 = c$ {\it a.s.}\\
Given (\ref{condeqn}), (\ref{eqn9})
holds (even if one $X_i, ~i=1, 2$ does not have finite variance).
}
\begin{remark}
So, if (\ref{eqn9}) holds, then it will hold only either as $\var{Y} = \var{X_1}$ or $\var{Y} = \var{X_2}$. Similarly, we can characterize the situation 
where $\var{Z}$ $= \var{X_1} + \var{X_2}$.
\end{remark}
Write $X(X_1,X_2) := \min(X_2, \max(X,X_1))$ for truncating a random variable $X$ by $X_2$
from above and by $X_1$ from below. Using Theorem \ref{var-bound1}, we 
have,
\begin{corollary}
If $X$, $X_1$ and $X_2$ have finite second moments, then 
\begin{equation}
\var{X(X_1,X_2)} \leq \var{X_2} + \var{X_1} + \var{X}
\end{equation}
\end{corollary}
As in Theorem \ref{var-eqt}, we would like to find conditions when there could 
be an equality in this bound.
\begin{theorem}
\label{var-eqt2}
Suppose $X$, $X_1$ and $X_2$ have finite variances.
\begin{equation}
\label{equality}
\var{X(X_1,X_2)} = \var{X_2} + \var{X_1} + \var{X}
\end{equation}

iff, for some constants $c_1$ and $c_2$, one of the following sets of conditions hold:
\begin{equation}
\label{eqn13}
\begin{array}{cl}
(1) & c_2 = X_2 \geq X \geq X_1 = c_1 \mbox{\it\ a.s.} \\
(2) & c_2 = X \geq X_1 = c_1  \mbox{ and } c_2 = X \geq X_2 \mbox{\it\ a.s.} \\
(3) & c_2 = X_2 \geq X_1 \geq X = c_1   \mbox{\it\ a.s.} \\
(4) & c_2 = X_1 \geq X = c_1 \mbox{ and } c_2 = X_1 \geq X_2 \mbox{\it\ a.s.}\\ 
\end{array}
\end{equation}
\end{theorem}
\Proof{
We will use Corollary \ref{cor1} in the following. Equation (\ref{equality}) holds, iff\\
\begin{eqnarray}
 \begin{array}{lccr}
\var{X(X_1,X_2)}& & & \\  
= \var{X_2} + \var{\max(X,X_1)} & \mbox{iff} & 
                                      \var{\min(X,X_1)} = 0  & \mbox{(I)}\\
  &            & \mbox{\it i.e. } \min(X,X_1) = c_1 \mbox{\it\ a.s.} & \\
  & \mbox{\it and} & X \geq X_1 = c_1 \mbox{\it\ a.s.} & \mbox{(a)} \\
  &            &      \mbox{\it or}          &            \\
  &            & X_1 \geq X = c_1 \mbox{\it\ a.s.} & \mbox{(b)} \\
  & & & \\
= \var{\min(X_2,\max(X,X_1))}&\mbox{iff}& 
                    \max(X_2, \max(X,X_1)) = c_2 \mbox{\it\ a.s.} & \mbox{(II)}\\
  & \mbox{\it and} & c_2 = X_2 \geq \max(X,X_1) \mbox{\it\ a.s.} & \mbox{(c)} \\
  &            &      \mbox{\it or}         &                       \\
  &            & X_2 \leq \max(X,X_1)=c_2 \mbox{\it\ a.s.} & \mbox{(d)} \\
\end{array} \nonumber
\end{eqnarray}
Subsuming the {\it and} in conditions (I) and (II) to the {\it or} of conditions (a,b) and (c,d) respectively, we get the four conditions 
(1,2,3 and 4) mentioned in the statement of the theorem.\\
If (\ref{eqn13}) holds, then (\ref{equality}) follows (even if some variances
are not finite).
}

\begin{remark}
Conditions (1) and (3) are `symmetric' w.r.t. roles of $X$ and $X_1$; so are the other two, the reason being $\max(X,X_1)$ in the definition of $X(X_1,X_2)$ is symmetric about its arguments. 
\end{remark}

Finally, let $s_2 > s_1$ where $s_1$ and $s_2$ are constants 
and consider $Y_1 := \min(X,s_1)$ and $Y_2 := \min(X,s_2)$.
Since $Y_1 = \min(Y_2, s_1)$, we have $ \var{Y_1} \leq \var{Y_2}$.
Thus $Y(\cdot)$ is a non-decreasing 
function of $s$. Similarly, $Z(s) := \max(X,s)$ is a non-increasing function of $s$; 
see also \cite{Chow69}.

\small
\section*{\small Acknowledgments} We thank Prof. Narayan Rangaraj of IITB for 
discussions on some inventory models that lead to this work. We also thank 
Prof. Vivek Borkar of TIFR for useful comments.


\normalsize

\begin{thebibliography}{9}
\bibitem{Asmussen}
Asmussen, S. {\em Applied probability and queues},
John Wiley \& Sons, Chichester, 1987.

\bibitem{Bertsekas}
Bertsekas, D. P. {\em Dynamic programming and optimal control, Vols. 1 and 2},
Athena Scientific, Belmont, 1995.

\bibitem{Chow69}
Chow, Y. S. and Studden, W. J., 
Monotonicity of the variance under truncation and variations of Jensen's 
inequality, {\it The Ann. of Math. Stat.}, {\bf 40}, 1106-1108, 1969.

\bibitem{Chow88}
Chow, Y. S. and Teicher, H., 
{\em Probability theory: Independence, interchangeability and martingales}, 
Springer-Verlag, New York, 2nd ed., 1988.


\bibitem{Graves} 
Graves, S. C., Rinnooy Kan, A. H. G., and Zipkin, P. H., (Editors) 
{\em Logistics of production and inventory, Handbooks in Operations Research 
and Management Science}, Vol. 4, North-Holland, Amsterdam, 1993.

\bibitem{Iglehart}
Iglehart, D. L., Optimality of {$(s, S)$} policies in the infinite horizon 
dynamic inventory problem. {\it Management Sci.}, {\bf 9}, 259-267, 1963. 

\bibitem{Lindley}
Lindley, D. V., On the theory of queues with a single server, {\it Proc. Camb. 
Philos. Soc.}, {\bf 48}, 277-289, 1952.  

\bibitem{Papadatos}
Papadatos, N., Expectation bounds on linear estimators from dependent samples, {\it Journal of Statistical Planning and Inference}, {\bf 93}, 17-27, 2001. 

\bibitem{Scarf}
Scarf, H. E., The optimality of {($s, S$)} policies in the dynamic inventory 
problem, {\em Mathematical methods in social sciences} Editors: K. A. Arrow, 
S. Karlin and P. Suppes, Stanford University Press, Stanford, 1960. 

\bibitem{Wolff}
Wolff, R. W.,
{\em Stochastic modeling and the theory of queues},
Prentice-Hall International Inc., Englewood Cliffs, 1989.

\end{thebibliography}
\end{document}